\documentclass[12pt,leqno]{article}
\usepackage{amsmath,amssymb,latexsym}
\usepackage{amscd}
\usepackage{a4}
\usepackage{mathrsfs}
\usepackage{srcltx}
\usepackage{theorem}
\usepackage{amsfonts}

\newcommand{\R}{\mathbb{R}}        
\newcommand{\C}{\mathbb{C}}        
\newcommand{\N}{\mathbb{N}}        
\newcommand{\Rn}{\mathbb{R}^d}     
\newcommand{\Cn}{\mathbb{C}^d}
\newcommand{\CN}{C^\infty(\Rn)}

\newcommand{\Proof}{\textbf{Proof:} \ }
\newcommand{\qed}{\hspace*{\fill} $\Box $}

\newcommand{\cA}{\mathscr{A}}
\newcommand{\cD}{\mathscr{D}}
\newcommand{\cE}{\mathscr{E}}
\newcommand{\cH}{\mathscr{H}}

\newcommand{\cC}{\mathscr{C}}
\newcommand{\cW}{\mathscr{W}}

\newcommand{\cS}{\mathscr{S}}

\newcommand{\EP}{\cE'(\Rn)}

\newcommand{\db}{\bar{\partial}^\be}

\newcommand{\be}{{\bf 1}}

\renewcommand{\span}{\mathrm{span}}

\newcommand{\om}{\omega}

\newcommand{\ssubset}{\subset\subset}

\newtheorem{proposition}{Proposition}[section]
\newtheorem{example}[proposition]{Example}
\newtheorem{lemma}[proposition]{Lemma}
\newtheorem{corollary}[proposition]{Corollary}
\newtheorem{theorem}[proposition]{Theorem}

\newtheorem{definition}{Definition}

\newcommand{\Ao}{\mbox{$\cA (\Omega)$}}

\newcommand{\ind}{\mbox{\rm ind}}

\newcommand{\Om}{{\Omega}}

\newcommand{\eps}{\varepsilon}

\newcommand{\conv}{{\rm conv}}
\newcommand{\mconv}{{\rm mconv}}

\newcommand{\Exp}{{\rm Exp}}
\newcommand{\Log}{{\rm Log}}

\newcommand{\supp}{\mathrm{supp}\, }
\newcommand{\lp}{\mbox{\rm lim proj}}

\newcommand{\vp}{\varphi}

\title{{\sc $\cE'$ as an algebra by multiplicative convolution
}}
\author{Dietmar Vogt}
\date{}

\begin{document}

\maketitle

\footnotetext{
    \hfil\break\indent \begin{minipage}[t]{14cm}{\em Key words and phrases:} Operators on smooth functions, distributions with compact support, algebra of operators, monomials as eigenvectors, Hadamard operators, Euler operators. \end{minipage}.
     \hfil\break\indent}
\centerline{\em Dedicated to the memory of Pawe{\l} Doma\'nski, great mathematician and friend}

\begin{abstract}We study the algebra $\cE'(\Rn)$ of distributions with compact support equipped with the multiplication $(T\star S)(f)=T_x(S_y(f(xy))$ where $xy=(x_1 y_1,\dots,x_d y_d)$. This allows us a very elegant access to the theory of Hadamard type operators on $C^\infty(\Omega)$, $\Omega$ open in $\Rn$, that is, of operators which admit all monomials as eigenvectors. We obtain a representation of the algebra of such operators as an algebra of holomorphic functions with classical Hadamard multiplication. Finally we study global solvability for such operators, in particular of Euler differential operators, on open subsets of $\R_+^d$.
\end{abstract}

In the present paper we continue the investigations in \cite{Vhad}. For that purpose we define the \em multiplicative convolution \rm in the space $\EP$ of distributions of compact support on $\Rn$ by $(T\star S)g =T_x S_y f(xy)$. Equipped with this multiplication $\EP$ is a commutative algebra. Each $T\in\EP$ defines a convolution operator $N_T:S\mapsto T\star S$ in $\EP$. For open $\Om,\,\Om'\subset\Rn$ we show that $N_T \cE'(\Om')\subset \cE'(\Om)$ if and only if $\supp T\cdot\Om'\subset \Om$, in particular, $N_T$ defines an operator in $L(\cE'(\Om))$ if and only if $\supp T\subset V(\Om)$ where $V(\Om)$ denotes the set of dilations acting on $\Om$. In this case $N_T=M_T^*$ where $M_T\in L(\cE(\Om))$ is the Hadamard operator on $\CN$ defined by $M_T f(y)=T_x f(xy)$ that is $N_T$ is the dual operator of $M_T$. This allows us an elegant access to the theory of Hadamard operators on $\cE'(\Om)$ for open sets $\Om\subset\Rn$, that is, of operators which allow all monomials as eigenvectors. These have been studied and characterized in \cite{Vhad}. The algebra $M(\Om)$ of such operators is a closed subalgebra of $L(C^\infty(\Om))$. It is isomorphic to the algebra $(\cE'(V(\Om)),\star)$ and we determine the topology induced from $L_b(C^\infty(\Om))$ on $\cE'(V(\Om))$. Then we use the exponential diffeomorphism between $Q_+:=(0,+\infty)^d$ and $\Rn$ to study, by transfer from the classical theory, global solvability of Hadamard operators on $C^\infty(\Om)$ where $\Om$ is an open subset of $Q_+$.  Hadamard-type operators assigned to distributions with support $\{(1,\dots,1)\}$ are differential operators of the form $P(\theta)$, where $P$ is a polynomial and $\theta_j=x_j\,\partial/\partial x_j$. They are called Euler operators and we study, in particular, global solvability for such operators on open subsets of $Q_+$. Euler operators on general open subsets of $\Rn$ are studied, by different methods, in Doma\'nski-Langenbruch \cite{DLeul}, global solvability in $\cS'(\Rn)$ is shown in \cite{V2}.  Finally we show that the algebra $M(\Om)\cong (\cE'(V(\Om)),\star)$ is isomorphic to an algebra of holomorphic functions around zero where the multiplication is the classical Hadamard multiplication, that is, multiplication of the Taylor coefficients.

Hadamard operators on the space $\Ao$ of real analytic functions have been studied and characterized in Doma\'nski-Langenbruch \cite{DLre}, \cite{DLalg}, \cite{DLhad}. Analogous problems as in the present paper for real analytic functions are studied in Doma\'nski-Langenbruch-Vogt \cite{DLV}. While the problems are analogous, the results, the methods and the difficulties to overcome are, in part, quite different.
For further references and also for the background of the problem we refer to \cite{Vhad}.

We use standard notation of Functional Analysis, in  particular, of distribution theory. For general information and unexplained notation we refer to \cite{DK}, \cite{MV}, \cite{LSI}. For the classical product of convolution in $\cE'(\Rn)$ see \cite{DK} and \cite{LS}.

\section{Multiplicative convolution on $\cE'$}

On $\cE'(\R^d)$ we define the \em multiplicative convolution \rm:

\begin{definition} For $T,\,S\in \cE'(\R^d)$ and $f\in C^\infty(\R^d)$ we set $(T\star S)f=T_x(S_y f(xy))$.
\end{definition}

We freely use this notation, in different context, also for other combinations of $T,\,S$ and $f$, whenever the right hand side makes sense. As multiplication on $\cE'(\R^d)$ it has the following basic properties:

\begin{proposition}\label{pr1} 1. $(\cE'(\R^d),\star)$ is a commutative algebra.\\
2. $\delta_\be$ is the unit element.\\
3. The multiplication is hypocontinuous.\\
4. $\supp (T\star S)\subset \supp T\cdot \supp S$.
\end{proposition}

\Proof 1. To show commutativity we remark that for any $\alpha\in\N_0^d$ we have
$$(T\star S)x^\alpha= Tx^\alpha\cdot Sx^\alpha$$
Therefore $T\star S$ and $S\star T$ coincide on the polynomials which are dense in $C^\infty(\R^d)$.

2. is trivial.
To show 3. we have, with suitable $K$ and $p$ for $T$, resp. $L$ and $q$ for $S$, and a constant $C>0$ depending only on $L,\,p,\,q$ the following estimates:
\begin{eqnarray*}
|(T\star S)f| &=&|T_x(S_y f(xy))|\\ &\le& \|T\|^\ast_{K,p} \sup_{x\in K,|\alpha|\le p}|S_y (y^\alpha f^{(\alpha)}(xy))|\\ &\le& \|T\|^\ast_{K,p}\,\|S\|^\ast_{L,q}\sup_{x\in K,|\alpha|\le p}\,\sup_{y\in L, |\beta|\le q}|\partial_y^\beta(y^\alpha(f^{(\alpha)}(yx))|\\ &\le& C\,\|T\|^\ast_{K,p}\,\|S\|^\ast_{L,q} \sup_{z\in K L,|\gamma|\le pq}|f^{(\gamma)}(z)|.
\end{eqnarray*}
This shows 3. and implies 4. \qed

Let $\Om'$ and $\Om$ be open subsets of $\R^d$. From Proposition \ref{pr1} we obtain immediately
\begin{corollary}\label{c1} If $T\in\cE'(\R^d)$ and $\supp T\cdot\Om'\subset\Om$, then $N_T:S\mapsto T\star S$ defines a map in $L(\cE'(\Om'),\cE'(\Om))$.
\end{corollary}

\Proof We have only to remark that $\cE'(\Om')$ is bornological and that, by Proposition \ref{pr1}, the map $N_T$ sends equicontinuous (=bounded) sets into equicontinuous (=bounded) sets. \qed

We want now to describe the maps $N\in L(\cE'(\Om'),\cE'(\Om))$ for which there is $T\in \cE'(\R^d)$ with $N=N_T$. We denote the set of these maps by $N(\Om',\Om)$. We set $\R_*=\R\setminus \{0\}$.

\begin{lemma}\label{l1} For $N\in L(\cE'(\Om'),\cE'(\Om))$ the following are equivalent: \\
 1. $N(S_1)\star S_2=S_1\star N(S_2)$ for all $S_1,\,S_2\in\cE'(\Om')$\\
 2. There is $T\in\cE'(\R^d)$ such that $N(S)=T\star S$ for all $S\in\cE'(\Om')$.
\end{lemma}

\Proof To show 1. $\Rightarrow$ 2. we choose $\eta\in\Om'\cap \R_*^d$. Then for all $S\in\cE'(\Om')$ we have $N(S)\star \delta_\eta= S\star N(\delta_\eta)$ and therefore $N(S) = T\star S$ with $T=N(\delta_\eta)\star \delta_{1/\eta}$. The reverse direction is obvious.\qed

For two sets $X,\,Y\subset \R^d$ we set $V(X,Y):=\{\eta\in\R^d\,:\,\eta X\subset Y\}$, $V(X)=V(X,X)$ and $Y_x=\{\eta\in \R^d\,:\,\eta x\in Y\}$.

\begin{lemma}\label{l2} If $X$ is open and $Y$ is closed, then $V(X,Y)$ is closed and $V(X,Y)=\bigcap_{x\in X\cap \R_*^d} \frac{1}{x}Y$. If, in particular, $Y$ is compact, then also $V(X,Y)$.
\end{lemma}

\Proof Because $X\cap \R_*^d$ is dense in $X$ we have
$$V(X,Y)=\{\eta\in\Rn\,:\,\eta\, (X\cap \R_*^d)\subset Y\}=\bigcap\nolimits_{x\in X\cap \R_*^d}\frac{1}{x} Y.$$
Since $\frac{1}{x} Y$ is closed for all $x\in \R_*^d$ the set $V(X,Y)$ is closed, same argument in the case of compact $Y$. \qed

\begin{theorem}\label{t1} For $T\in\cE'(\Rn)$ the following are equivalent:\\
1. $N_T\in L(\cE'(\Om'),\cE'(\Om))$.\\
2. $N_T\in L(C'(\Om'),\cE'(\Om))$.\\
3. $\supp T\subset V(\Om',\Om)$.
\end{theorem}

\Proof Since 3. $\Rightarrow$ 1. and 1. $\Rightarrow$ 2. are obvious we have only to show 2. $\Rightarrow$ 3.

If $N_T\in L(C'(\Om'),\cE'(\Om))$ we have for every $\eta\in\Om'$ that $\supp (T\star\delta_\eta)\subset \Om$.  We fix an open $\om'\ssubset\Om'$. Because of the continuity of $N_T$ there is a compact subset $L\subset\Om$ such that $\supp (T\star\delta_\eta)\subset L$ for all $\eta\in\om'$. For $\eta\in \R_*^d$ we have $\supp (T\star\delta_\eta)=\eta\, \supp T$. So we have shown that $\supp T\subset \bigcap_{\eta\in \om'\cap \R_*^d} \frac{1}{\eta}L=V(\om',L)\subset V(\om',\Om)$. The equality comes from Lemma \ref{l2}. Since $V(\Om',\Om)=\bigcap_{\om'\ssubset\Om'}V(\om',\Om)$ the proof is complete. \qed

\section{Transposed description}

Let $N_T\in N(\Om',\Om)$, then the transposed operator $N_T^\ast \in L(C^\infty(\Om),C^\infty(\Om'))$ is given by
$$(N_T^\ast f)(x) =\langle \delta_x, N_T^\ast f\rangle =\langle T\star \delta_x, f\rangle  = T_yf(xy)=:(M_T f)(x),$$ that is $N_T^\ast = M_T$. We extend this to a complete characterization:

\begin{theorem}\label{t2} For $M\in L(C^\infty(\Om),C^\infty(\Om'))$ the following are equivalent:\\
1. There is $T\in \cE'(V(\Om',\Om))$ such that $M=N_T^\ast$.\\
2. There is $T\in \cE'(V(\Om',\Om))$ such that $M=M_T$.\\
3. $M$ admits all monomials as eigenvectors.
\end{theorem}

\Proof The equivalence of 1. and 2. is obvious. Likewise 2. $\Rightarrow$ 3. is clear. We show 3. $\Rightarrow$ 1. We set $N:=M^\ast \in L(\cE'(\Om'),\cE'(\Om))$ and have to show that $N(S_1)\ast S_2 =S_1\star N(S_2)$ for all $S_1,\,S_2\in\cE'(\Om')$. Due to the density in $C^\infty(\Om')$ of the polynomials we have to show the equality only on the set of all monomials. For $f(x)=x^\alpha$ we obtain
$$(N(S_1)\star S_2)f=N(S_1)_{x}(S_2)_y (xy)^\alpha = S_1(M x^\alpha)S_2(x^\alpha)=m_\alpha S_1(x^\alpha) S_2(x^\alpha)$$
where $M x^\alpha = m_\alpha x^\alpha$. For $(S_1\star N(S_2))f$  we obtain the same, which shows the result. \qed

The operators described in Theorem \ref{t2} are called \em operators of Hadamard type, \rm or also \em multipliers \rm and they have been subject of many investigations where the equivalence of 2. and 3. was one of the main results.

\begin{definition} For open $\Om\subset\Rn$ we set  $M(\Om):=\{M\in L(C^\infty(\Om))\,:\, M \text{ Hadamard}\}.$
\end{definition}

We state two basic facts:
\begin{enumerate}
\item $M(\Om)$ is a closed subalgebra of $L_b(C^\infty(\Om)).$
\item  $T\mapsto M_T$ an algebra isomorphism $(\cE'(V(\Om)),\star)\to M(\Om)$.
\end{enumerate}
The first follows from $Mx^\alpha\subset\span\{x^\alpha\}$ for all $\alpha$, the second is a simple calculation.

 \section{Topologies} \label{Top}

 In this section we will study the problem: Which topology does the algebraic isomorphism $T\mapsto N_T$ from $\cE'(V(\Om))$ to $N_b(\Om',\Om)$ induce on $\cE'(V(\Om))$? Here $*_b$ denotes the topology of uniform convergence on bounded sets
 \bigskip

 The solution will be given in several steps. First we need some notation:

 For $\Om'\ssubset \Om$ open we set
 $$U=U_{\Om'}=V(\overline{\Om'},\Om)=\{\eta\in\Rn\,:\,\eta \Om'\ssubset\Om\}$$
and
$$NC_b(\Om',\Om)=\{N_T\,:\,T\in\cE'(U)\}\cap L_b(C'(\overline{\Om'}),\cE'(\Om)).$$
It is easily seen that $U_{\Om'}$ is an open neighborhood of $V(\Om)$.

\begin{proposition}\label{p3}
 $\cE'(U)\cong N_b(\Om',\Om)=NC_b(\Om',\Om)$ are topologically isomorphic. The first isomorphism is given by $T\mapsto N_T$.
 \end{proposition}

 \Proof From Proposition \ref{pr1} and its proof it is clear that $T\mapsto N_T$ sends every bounded set of $\cE'(U)$ into a bounded set of $L_b(\cE'(\Om'),\cE'(\Om))$. Since $\cE'(U)$ is bornological the map is continuous. Clearly the identical map $N_b(\Om',\Om)\to NC_b(\Om',\Om)$ is continuous.

Let $B\subset L_b(C'(\overline{\Om'}),\cE'(\Om))$ be bounded, then $\widetilde{B}=\{\vp(\delta_\eta)\,:\,\eta\in\overline{\Om'},\,\vp\in B\}$ is a bounded set in $\cE'(\Om)$. Therefore there is a compact set $K\subset\Om$ such that $\widetilde{B}$ is bounded in $\cE'(K)$.

 We fix now $\eta\in\overline{\Om'}\cap\R_*^d$ and obtain, by use of the proof of Lemma \ref{l1},
$$\{T\,:\,N_T\in B\}=\{\delta_{1/\eta}\star N_T(\delta_\eta)\,:\,N_T\in B\}\subset \{\delta_{1/\eta}\star S\,:\, S\in\widetilde{B}\}.$$
We set $L=\bigcap_{\eta\in\overline{\Om'}\cap\R_*^d}\frac{1}{\eta}K=V(\overline{\Om'},K)$ (see Lemma \ref{l2}). Then $L$ is a compact subset of $U$ and $\{T\,:\,N_T\in B\}$ is a bounded subset of $\cE'(L)$, since it is a bounded subset of $\cE'(\Rn)$ and the support of all of its elements is contained in $L$. So it is a bounded subset of $\cE'(U)$.

We have shown that the map $\mathcal{N}:\cE'(U)\to NC_b(\Om',\Om)$ is continuous and surjective. Moreover $\mathcal{N}^{-1}(B)$ is bounded for every bounded $B\subset L_b(C'(\overline{\Om'}),\cE'(\Om))$.
By Baernstein's Lemma (see \cite{MV}) $T\mapsto N_T$ is a topological imbedding of $\cE'(U)$ into $L_b(\cE'(\Om'),\cE'(\Om))$, which proves the result. \qed

Let now $\om_1\ssubset\om_2\ssubset\dots$ be an exhaustion of $\Om$.
We put $U_n:=U_{\om_n}$.
Then $U_1\supset U_2\supset\dots$ is a decreasing sequence of open neighborhoods of $V(\Om)$ with $\bigcap_n U_n= V(\Om)$.

\begin{definition}
 $\cE'_t(V(\Om)):=\lp_n \cE'(U_n).$
\end{definition}
Clearly the $t$-topology does not depend on the choice of the exhaustion.
$\cE'(U_n)$ is a (DF)-space hence  $\cE'_t(V(\Om))$ is a (PDF)-space.

We set  $N_b(\Om):= N_b(\Om,\Om)$ and $NC_b(\Om)=NC_b(\Om,\Om)$.

\begin{theorem}\label{t8}  $\cE'_t(V(\Om))\cong N_b(\Om)=NC_b(\Om)$ topologically, the first isomorphism is established by $T\mapsto M_T$.
\end{theorem}

\Proof  This follows immediately from Proposition \ref{p3} by forming projective limits. \qed

\section{Examples}
For the following examples see also \cite{Vhad}, for a more systematic treatment \cite{DLV}.
 We refer to the notation of Section \ref{Top}.

\begin{example} Let $\Om=\{(x,y)\in \R^2\,:\,1<y<2\}$ then $V(\Om)=\{(x,1)\,:\,x\in\R\}$.\\
 With $r_n\searrow 1$ set
$\om_n=\{(x,y)\in\R^2\,:\,|x|<n, r_n< y <2/r_n\}$. We obtain $U_n=\{(x,y)\,:\,1/r_n < y <r_n\}$.
\end{example}

\begin{example} Let $\Om=\{(x,y)\in \R^2\,:\,0<x,\,1<y<2\}$ then $V(\Om)=\{(x,1)\,:\,x>0\}$.
 With $r_n\searrow 1$ set
$\om_n=\{(x,y)\in\R^2\,:\,1/n<|x|<n, r_n< y <2/r_n\}$. We obtain $U_n=\{(x,y)\,:\,0<x,\,1/r_n < y <r_n\}$.
\end{example}

 In both cases $V(\Om)$ is closed in some $U_n$. If this is the case,
\begin{itemize}
 \item then it is closed  in all $U_m$ for $m\ge n$
 \item $\cE_t'(V(\Om))\cong \cE'(V(\Om))$ equipped with the canonical (DF)-topology inherited from $\cE'(U_n)$.
 \item $V(\Om)$ is locally compact, $\sigma$-compact and $\cE_t'(V(\Om))=\ind_{K\subset V(\Om)}\cE'(K)$.
 \end{itemize}

\medskip

We have shown:

\begin{theorem}\label{t9} If $V(\Om)$ is closed in some $U_n$, then $V(\Om)$ is locally compact, $\sigma$-compact and $M(\Om)\cong \cE'(V(\Om))$ with  its standard (DF)-topology.
\end{theorem}

\begin{example}
If $\overline{V(\Om)}\setminus V(\Om)$ is finite then the theorem applies.
\end{example}

We need the following property of dilation sets, see \cite{Vhad} or \cite[Proposition 4.4]{DLV}. We give the proof for the convenience of the reader.

 \begin{lemma}\label{l5}$y\in\overline{V(\Om)}\setminus V(\Om)$ then there is $j$ with $y_j=0$.
\end{lemma}

\Proof If $y\in \overline{V(\Om)}$ and $\min_j|y_j|\ge 2\eps>0$ then $d(yx,\partial\Om)\ge\eps d(x,\partial\Om)$ for all $x\in\Om$, since that holds for an approximating sequence in $V(\Om)$. Hence $yx\in\Om$. Therefore for $y\in\overline{V(\Om)}\setminus V(\Om)$ there is $j$ such that $y_j=0$. \qed

\begin{theorem}\label{t10} If $d=1$ then $V(\Om)$ is locally compact, $\sigma$-compact and $M(\Om)\cong \cE'(V(\Om))$ with  its standard (DF)-topology.
\end{theorem}

\Proof By Lemma \ref{l5} we have $\overline{V(\Om)}\setminus V(\Om)\subset\{0\}$. By Theorem \ref{t9} the result follows. \qed

\section{Relation to classical convolution}\label{s5}

We set $Q_+=\{x\in\R^d\,:\,x_j>0\text{ for }j=1,\dots,d\}$ and restrict our attention to open sets $\Om\subset Q_+$. The map $\Log(x)=(\log x_1,\dots,\log x_d)$ defines a diffeomorphism from $Q_+$ onto $\R^d$ whose inverse is $\Exp(x)=(\exp x_1,\dots,\exp x_d)$.

The map $C_\Log:\vp\to \vp\circ \Log$ is an isomorphism from $C^\infty(\Rn)$ onto $C^\infty(Q_+)$ whose inverse is $C_\Exp$, defined analogously. Then the transposed maps $C_\Log^\ast:\cE'(Q_+)\to\cE'(\Rn)$ and $C_\Exp^\ast:\cE'(\Rn)\to\cE'(Q_+)$ are isomorphisms.

The crucial observation relating the two kinds of convolution is:
\begin{lemma}\label{l3} $T\star S=C_\Exp^\ast((C_\Log^\ast T)\ast(C_\Log^\ast S))$ for all $T,\,S\in\cE'(Q_+)$.
\end{lemma}

A subset $X\subset Q_+$ will be called multiplicatively convex (mconvex) if with $x,\,y\in X$ also $x^t y^{1-t}\in X$ for all $0<t<1$.
For a subset $X\subset Q_+$ we define the multiplicatively convex hull $\mconv X$ as the smallest mconvex set, which contains $X$. We have $\mconv X = \Exp (\conv (\Log X))$ where $\conv$ denotes the convex hull. Explicitly we have
$$\mconv X =\Big\{\prod_j x_j^{\lambda_j}\,:\,x_j\in X, 0<\lambda_j, \sum_j\lambda_j=1\Big\}$$
where the sums and products are finite and the powers are taken coordinatewise.

We obtain the multiplicative analogue of the Theorem of Lions \cite{Lions}.
\begin{theorem}\label{t3} $\mconv(\supp (T\star S))=\mconv(\supp T)\cdot\mconv(\supp S)$ for all $T,\,S\in \cE'(Q_+)$.
\end{theorem}

For the following we might exploit Theorem \ref{t3}, but we prefer to use the direct transfer of properties of classical convolution operators.

First we consider the case of $\supp T=\{\be\}$. The operators $M_T$, resp. $N_T$, are then called \em Euler operators \rm. In this case $T$ has the form $T=\sum_{|\alpha|\le p} c_\alpha \delta_\be^{(\alpha)}$. Since $M_{\delta_\be^{(\alpha)}}=(-1)^{|\alpha|}x^\alpha \partial^\alpha$ the Euler operator $M_T$ has the form
$M_T=\sum_{|\alpha|\le p}b_\alpha x^\alpha \partial^\alpha.$

With $\theta_j = x_j \partial_j$ and different coefficients it can also be written as $M_T=\sum_{|\alpha|\le p} c_\alpha\, \theta^\alpha.$

Let now $\Om\subset Q_+$ be open. If we set, keeping the notation, $P(z)=\sum_{|\alpha|\le p} c_\alpha z^\alpha$, then $P(\partial) C^\infty(\Log\,\Om) = C^\infty(\Log\,\Om)$ if, and only if, $\Log\,\Om$ is $P(\partial)$-convex. This means that for every compact $K\subset \Log\,\Om$ there is a compact $L\subset\Log\,\Om$ with the following property: if $\vp\in\cD(\Log\,\Om)$ and $P(-\partial)\vp\in\cD(K)$ then $\vp\in\cD(L)$ or, equivalently: if $S\in\cE'(\Log\,\Om)$ and $P(-\partial)S\in\cE'(K)$ then $S\in\cE'(L)$ (see \cite{H}).

Now, it is obvious that this condition simply carries over to $P(\theta)$ and $\Om$, and we obtain:
\begin{theorem}\label{t4} $P(\theta) C^\infty(\Om)=C^\infty(\Om)$ if, and only if, $\Om$ is $P(\theta)$-convex, that is, for every compact $K\subset \Om$ there is a compact $L\subset\Om$ with the following property: if $\vp\in\cD(\Om)$ and $P(-\theta)\vp\in\cD(K)$ then $\vp\in\cD(L)$ or, equivalently: if $S\in\cE'(\Om)$ and $P(-\theta)S\in\cE'(K)$ then $S\in\cE'(L)$.
\end{theorem}

\Proof The first condition comes by transfer via $C_\Log$. The second condition implies the first one. To show its necessity we define the multiplication operator $M T= x_1\cdot..\cdot x_d\, T$.  Let $\theta^*$ denote the transpose of $\theta$, hence $P(\theta^*)$ the transpose of $P(\theta)$. From elementary calculations we get $M^{-1} P(-\theta) M = P(\theta^*)$. Since $M$ preserves supports on $\cE'(Q)$ the assertion follows from functional analytic principles. \qed

Using Theorem \ref{t3} we get:
\begin{corollary}\label{c2} If $\Om$ is mconvex then all Euler operators are surjective in $C^\infty(\Om)$.
\end{corollary}

Staying with mconvex $\Om\subset Q_+$ we study a more general situation. We fix $T\in\cE'(Q_+)$ and set $\Om'=V(\supp T,\Om)=V(\mconv\,\supp T,\Om)$. It is easy to see that $\Om'$ is open and mconvex. Clearly $\supp T\subset V(\Om',\Om)$, that is, $N_T\in L(\cE'(\Om'),\cE'(\Om))$.
Let $K\subset \Om$ be compact and mconvex. Assume moreover that $S\in\cE'(Q_+)$ and $\supp(T\star S)\subset K$, then also
$\mconv\,\supp T\cdot \mconv\,\supp S=\mconv\,\supp(T\star S)\subset K$. This implies that $\supp S\subset L:= V(\mconv\,\supp T, K)$.
Since $\supp T\subset \R_*^d$ we have $L\subset \frac{1}{\eta}´K$ for any $\eta\in\supp T$. Therefore $L$ is compact. We have shown:

\begin{lemma}\label{l4} If $\Om\subset Q_+$ is open and mconvex, $T\in \cE'(Q_+)$ and $\Om'=V(\supp T, \Om)$, then for every mconvex
compact set $K\subset\Om$ there is a compact set $L\subset \Om'$ such that the following holds: if $S\in\cE'(Q_+)$ and
$\supp (T\star S)\subset K$ the $\supp S\subset L$.
\end{lemma}

To get solvability conditions for our multplicative convolution equations we need the equivalent of an elementary solution. We set
for $T\in \cE'(Q_+)$
$$\breve{T}(z):= T_x(x^{-iz}),\quad z\in\C^d$$
and remark that $\breve{T}=\widehat{C_\Log^\ast T},\,\, \widehat{ }$\, denoting  the Fourier-transform.

\begin{definition} An entire function $J$  is said to be slowly decreasing if there exist positive
numbers $a, b, c$ such that for each point $x\in\Rn$ we can find a point $y\in\Rn$ with\\
(1) $|x-y|\le a \log(1 + |x|),$\\
(2) $|J(y)| \ge b/(1 + |y|^c).$
\end{definition}

This notation is due to Ehrenpreis and we obtain the following multplicative analogue to the Theorem of Ehrenpreis on completely inversible operators \cite{LE}. Note that in 2. below $(T\ast E)\vp = T_x(E_y \vp(xy))$ makes sense since for $x\in Q_+$ and $\vp\in\cD(Q_+)$ the function $x\mapsto \vp(xy)$ is in $\cD(Q_+)$.
\begin{theorem}\label{t5} For $T\in\cE'(Q_+)$ the following are equivalent\\
1. $M_T$ is surjective in $C^\infty(Q_+)$.\\
2. There exists $E\in \cD'(Q_+)$ such that $T\star E=\delta_\be$.\\
3. $\breve{T}$ is slowly decreasing.
\end{theorem}

\Proof 1. is equivalent to the surjectivity of $f\mapsto (C_\Log^\ast T)\ast f$ in $C^\infty(\Rn)$. By the Theorem of Ehrenpreis this is equivalent to $\widehat{C_\Log^\ast T}$ being slowly decreasing an this is equivalent to 3.

3. means that $\widehat{C_\Log^\ast T}$ is slowly decreasing and this is, by the Theorem of Ehrenpreis, equivalent to the existence of a distribution $W\in\cD'(\Rn)$ with $(C_\Log^\ast T)\ast W=\delta$. Given 3. we set $E=C_\Exp^\ast  W$ and obtain
$$T\star E=C_\Exp^\star((C_\Log^\ast T)\ast(C_\Log^\star E))=C_\Exp^\star((C_\Log^\ast T)\ast W)=C_\Exp^\star \delta=\delta_\be.$$
Given 2. we set $W=C_\Log^\ast E$ and proceed like before. \qed

Putting things together we obtain
\begin{theorem}\label{t6} If $\Om\subset Q_+$ is mconvex and open, $T\in\cE'(Q_+)$, $\breve{T}$ slowly decreasing and $\Om'=V(\supp T,\Om)$, then $M_T:C^\infty(\Om')\to C^\infty(\Om)$ is surjective.
\end{theorem}

\Proof First we remark that, by Theorem \ref{t5} $M_T \cE(\Om')$ is dense in $\cE(\Om)$. Due to the Surjectivity Criterion we have to show that for any bounded set $B\subset \cE'(\Om)$ the set $N_T^{-1}B$ is bounded in $\cE'(\Om)$. There is compact, mconvex $K\subset\Om$ such that $B\subset\cE'(K)$ and bounded there. By Lemma \ref{l4} there is compact $L\subset \Om'$ such that $N_T^{-1}\cE'(K)\subset \cE'(L)$. Then $\star$-convolution with $E$ (notation of Theorem \ref{t5}, 2.) yields a continuous linear map $R(N_T)\cap\cE'(K)\to \cE'(L)$ inverting $N_T$. So $N_T^{-1}B = E\star (R(N_T)\cap B)$ is bounded in $\cE'(L)\subset\cE'(\Om')$. \qed

An interesting fact is the following: for fixed $z\in\C^d$ we have $M_T (\xi\to\xi^z)[x]=T_y (xy)^z=T_\xi(\xi^z)x^z$. So all functions $x^z$ are eigenfunctions of the operator $M_T$ with eigenvalue $\breve{T}(iz)$. In particular we have $m_\alpha=\breve{T}(i\alpha)$ where $m_\alpha, \,\alpha\in\N_0^d$ is the multiplier sequence for the Hadamard operator $M_T$.

This all can be done for each of the $2^d$ quadrants, but not for the whole euclidean space, as the following simple example shows. Take the operator $x\frac{\partial}{\partial x}$. It is easily seen to be surjective on $C^\infty(Q_{+/-})$, but it is not surjectice on $C^\infty(\R)$ since $xf'(x)$ is $0$ in $0$ for any $f$.

\section{Laurent representation theorem}\label{s6}

For $T\in\cE'(\Rn)$ and $z\in\Cn,\,z_j\neq 0$ for all $j$ we set
$\cC(z)=\prod_{j=1}^d\frac{1}{z_j}.$ For any subset $B\subset \Rn$ we define
$$\cW(B)=\{z\in\Cn\,:\,\xi_j\neq z_j\text{ for all }\xi\in B\text{ and }j=1,\dots,d\}.$$

Since $\cC$ defines a distribution on $\Cn$ also $\cC_T :=T\ast \cC$ is a distribution on $\Cn$. For $z\in\cW(\supp T)$ we obtain $\cC_T(z)=(T\ast \cC)(z)=T_\xi(\cC(z-\xi))$ which is a holomorphic function on $\cW(\supp T)$.

On $\cD'(\Cn)$ we set $\bar{\partial}^\be:=\pi^{-d} \frac{\partial}{\partial \bar{z}_1}\cdots\frac{\partial}{\partial \bar{z}_d}$, then $\cC$ is a fundamental solution for $\db$ and therefore $\db\cC_T=T$.

If $\supp T\subset \{\xi\in\Rn\,:\,|\xi|_\infty\le R\}$ then $(\C\setminus[-R,+R])^d\subset \cW(\supp T)$ and $\cC_T$ extends to a holomorphic function on $(\widehat{\C}\setminus[-R,+R])^d$, $\widehat{\C}$ denoting the Riemann sphere. $\cC_T(z)=0$ outside $\Cn$. $\cC_T$ is the unique solution of $\db u=T$ with this property.

For $\min_j |z_j|> R$ the function $\cC_T(z)$ is defined and holomorphic and it has the expansion
$$\cC_T(z)=\frac{1}{z_1\cdots z_d}\sum_{\alpha\in\N_0^d} T_\xi(\xi^\alpha) \frac{1}{z^\alpha}=\sum_{\alpha\in\N_0^d} m_\alpha \frac{1}{z^{\alpha+\be}}.$$

We have proved the following
\begin{proposition}\label{p2} Let $B\subset\Rn$ be compact and closed with respect to multiplication. The algebra $(\cE'(B),\star)$ is algebra-isomorphic to the algebra of all distributions on $\Cn$ which are holomorphic on $(\C\setminus[-R,+R])^d$ for some $R>0$, regular with value 0 in all infinite points of $\widehat{\C}^d$ and zero solutions for $\db$ outside $B$, equipped with Hadamard multiplication of the coefficients of their Laurent expansion around $(\infty,\dots,\infty)$.
\end{proposition}

We will now give a description in terms of properties of the functions on $(\C\setminus\R)^d$. All of the following is well known (see \cite{Sch,Ti}). We indicate the proofs for the convenience of the reader.

\begin{theorem}\label{t7}
1. If $f=\cC_T$ on $(\C\setminus\R)^d$ for $T\in\cE'(\Rn)$ then there is $p\in\N_0$ such that
\begin{equation}\label{e1}\sup_{z\in(\C\setminus\R)^d} |f(z)|\,|y_1\cdots y_d|^p<+\infty\end{equation}
where $z=x+iy$.

2. If $f\in H((\C\setminus [-R,+R])^d)$ for some $R>0$ is regular with value 0 in all infinite points of $\widehat{\C}^d$ and fulfills (\ref{e1}), then there is $T\in  \cE'(\Rn)$ such that $f=\cC_T$.

3. If $f=\cC_T$ then for $x\in\Rn$ we have: $x\not\in\supp T$ if, and only if there is a complex neighborhood $\om$ of $x$ such that $f=f_1+\dots +f_d$ on $\om\cap (\C\setminus\R)^d$ and $f_j \in H(\om_j)$ where $\om_j=\om\cap\{z=x+iy\in\Cn\,:\, y_\nu\neq 0 \text{ for }\nu\neq j\}$.
\end{theorem}

\Proof 1. follows directly from the continuity estimates.

2. The distribution $T$ is given by
\begin{equation}\label{e2}T\vp=\lim_{y\to 0+}\sum_{e}\mathrm{sgn}(e) \int \vp(x) f(x+iey) dx \end{equation}
for $\vp\in\cD(\Rn)$. Here $e\in\{+1,-1\}^d$, $\mathrm{sgn}(e)=e_1\cdots e_d$ and $y\to 0+$ means $y_j\to 0+$ for all $j$.

3. Sufficiency of the condition follows directly from (\ref{e2}). Necessity follows from the fact that solutions of $\db u=0$ have locally the form $u=u_1+\dots+u_d$ where $\frac{\partial u_j}{d\bar{z}_j}=0$. \qed

We collect all this information in a theorem.

\begin{theorem}\label{t11}
$T\mapsto \cC_T$ defines an algebra isomorphism of $M(\Om)\cong \cE'(V(\Om))$ to the following algebra $\cH_C(\Om)$ of holomorphic functions with Hadamard multiplication of Laurent coefficients around $(\infty,\dots,\infty)$
\begin{enumerate}
\item $f\in H((\widehat{\C}\setminus [-R,+R])^d)$ for some $R>0$, $f$ is regular with value 0 in all points outside $\Cn$.
\item There is $p\in\N_0$ such that
\begin{equation*}\sup_{z\in(\C\setminus\R)^d} |f(z)|\,|y_1\cdots y_d|^p<+\infty\end{equation*}
where $z=x+iy$.
\item $\sigma (f)$ is compact and $\sigma (f)\subset V(\Om)$.
\end{enumerate}
Here $\sigma(f)$ denotes the hyperfunction support of $f$.
\end{theorem}

We recall the definition:

\begin{definition} \label{d1}$x\not\in \sigma(f)$ if in a complex neighborhood of $x$ we have $f=f_1+\dots+f_d$ where $f_j$ extends in the $j$-th variable holomorphically across $\R$.
\end{definition}

\section{Hadamard representation theorem}

In a next step we want to change the equivalence into one with Hadamard multiplication of power series. For that we observe that for $z\in \C_*^d\cap\cW(\supp T)$ we have
$$\cC_T(z)= \frac{1}{z_1\cdots z_d}\; T_\xi\Big(\prod_j\frac{1}{1-\xi_j\cdot1/z_j}\Big)$$
and therefore for $z\in \C_*^d$ such that $1/z\in\cW(\supp T)$
$$T_\xi\Big(\prod_j\frac{1}{1-\xi_j z_j}\Big)=\frac{1}{z_1\cdots z_d}\; \cC_T(1/z).$$

\bigskip

We set $$C_T(z)=T_\xi\Big(\prod_j\frac{1}{1-\xi_j z_j}\Big)$$
on
$$W(\supp T) :=\{z\in\Cn\,:\,x_j z_j\neq 1\text{ for all }x\in \supp T\text{ and }j=1,\dots,d\}.$$
We have $(\C\setminus\R\cup [-R,+R])^d\subset W(\supp T)$ for some $R>0$ and for $|z|_\infty<R$
$$C_T(z)=\sum_\alpha m_\alpha x^\alpha.$$

This establishes an algebra isomorphism of $M(\Om)\cong \cE'(V(\Om))$ to an algebra of holomorphic functions with Hadamard multiplication, whose properties we get from the previous case by the formula
$$C_T(z)=\frac{1}{z_1\cdots z_d}\,\cC_T(1/z).$$

\vspace{.5cm}

\noindent Bergische Universit\"{a}t Wuppertal,
\newline Dept. of Math., Gau\ss -Str. 20,
\newline D-42119 Wuppertal, Germany
\newline e-mail: dvogt@math.uni-wuppertal.de

\end{document}